\documentclass[12pt]{article}
\usepackage{graphics}
\usepackage{amsmath}
\usepackage{amssymb}
\usepackage{graphicx}
\usepackage{makeidx}
\usepackage{mathrsfs}
\usepackage{amsfonts}
\usepackage[mathscr]{eucal}
\usepackage{amsmath}

\DeclareGraphicsExtensions{.eps,.pdf,.jpg,.png}

\def\RR{\vbox {\hbox to 8.9pt {I\hskip-2.1pt R\hfil}}}

\def\pni{\par\noindent}
\def\vsh{\smallskip}

\def\vsp{\vsh\pni} %% ie. \smallskip + \par

\numberwithin{equation}{section}

\begin{document}
%%%%%%%%%%%%%%%%%%%%%%%%%%%%%%%%%%%%%%%%%%%%%%%%%%%%%%%%%%%%%%%%%%%%%%%%
\font\title=cmbx12 scaled\magstep2
\font\bfs=cmbx12 scaled\magstep1
\font\little=cmr10
\begin{center}
{\title Transient waves in linear dispersive media  \\
with dissipation: an approach based\\
  on the  steepest descent path}
 \\  [0.50truecm]
 Francesco MAINARDI$^{(1)}$,
 Andrea MENTRELLI$^{(2), (3)}$ and
 \\
 Juan Luis GONZ{\'A}LEZ  SANTANDER$^{(4)}$
\\[0.25truecm]
$^{(1)}$  
 Department of Physics and Astronomy, University of Bologna and INFN.
\\ {Via Irnerio 46, I-40126 Bologna, Italy};
\\{ francesco.mainardi@unibo.it; mainardi@bo.infn.it; fracalmo@gmail.com}
\\ 
$^{(2)}$ Department of Mathematics and AM$^2$,  University of Bologna. 
\\ Via Saragozza 8, 40123 Bologna, Italy;
 {andrea.mentrelli@unibo.it}
\\
$^{(3)}$ Istituto Nazionale di Fisica Nucleare (I.N.F.N.), Sezione di Bologna,
\\  I.S. FLAG,  Viale Berti Pichat 6/2, 40127 Bologna, Italy;
\\
$^{(4)}$  Department of Mathematics, University of  Oviedo.
 \\ C Leopolodo Calvo Sotelo 18, 33007 Oviedo, Spain;
 {gonzalezmarjuan@uniovi.es}
%%
%\vskip 0.25truecm %% [0.25truecm]
%{\bf Version of \today}
 \vskip 0.25truecm
 % [0.25truecm]
 {\bf Published in Mathematics (MDPI) 
  Vol.13 No 21 (2025), 3418}
  \\ {\bf DOI: 
10.3390/math13213418}
\end{center}
%%%%%%%%%%%%
\begin{abstract}
In the study of linear dispersive media it is of primary interest to gain knowledge of the impulse response of the material. %, from which the response to arbitrary inputs may be obtained by means of the convolution operator.
	The standard approach to compute the response involves a Laplace transform inversion, i.e., the solution of a Bromwich integral, which can be a notoriously troublesome problem. % when analytical techniques fall short and a numerical solution is the only way to go.
	In this paper we propose a novel approach to the calculation of the impulse response, based on the well assessed method of the steepest descent path, which results in the replacement of the Bromwich integral with a real line integral along the steepest descent path.
	In this exploratory investigation, the method is explained and applied to the case study of the Klein-Gordon equation with dissipation, for which analytical solutions of the Bromwich integral are available, as to compare the numerical solutions obtained by the newly proposed method to exact ones.
	Since the newly proposed method, at its core, consists in replacing a Laplace transform inverse with a potentially much less demanding real line integral, the method presented here could be of general interest in the study of linear dispersive waves in presence of dissipation, as well as in other fields in which Laplace transform inversion come into play.

\newpage

 \vsp {\bf Keywords}:
{Transient waves in linear viscoelasticity; Klein-Gordon equation with dissipation; Laplace transform; steepest descent method.}

\vsp {\bf Mathematics Subject Classification (MSC)}:
{41A60, 30E15, 44A10, 35L20, 33C10}

\end{abstract}

%
% \maketitle

%%%%%%%%%%%%%%%%%%%%%%%%%%%%%
%\begin{document}

%%%%%%%%%%%%%%%%%%%%%%%%%%%%%%%%%%%%%%%%%%
%\setcounter{section}{-1} %% Remove this when starting to work on the template.
%\section{How to Use this Template}

\section{Introduction}

For uniaxial waves in an initially quiescent semi-infinite medium ($x\geq 0$%
), the response $r\left( x,t\right) $ of the medium to a pulse $r_{0}\left(
t\right) =r\left( 0,t\right) $ is \cite{MainardiKlein}:%
\begin{equation}
r\left( x,t\right) =\frac{1}{2\pi i}\int_{Br}\tilde{r}_{0}\left( s\right)
\exp \left( s\left[ t-n\left( s\right) \,x/c\right] \right) \,ds,
\label{r(x,t)_inv_Laplace}
\end{equation}%
where $s$ is the Laplace parameter, $\tilde{r}_{0}\left( s\right) $ is the
Laplace transform of the pulse $r_{0}\left( t\right) $, $n\left( s\right) $
is the Laplace transform of the medium refraction index, $Br$ is the
Bromwich path, and $c $ is the wave front velocity. When $r_{0}\left(
t\right) =\delta \left( t\right) $, its Laplace transform reads $\tilde{r}%
_{0}\left( s\right) =1$ so the corresponding response (the propagator or the
impulse response) is obtained inverting the Laplace transform:%
\begin{equation}
\tilde{r}_{\delta }\left( x,s\right) =\exp \left( -s\,n\left( s\right)
\,x/c\right) .  \label{r_delta_laplace}
\end{equation}%
Then, the solution corresponding to the generic $r_{0}\left( t\right) $ is
obtained by convolution:%
\begin{equation}
r\left( x,t\right) =r_{0}\left( t\right) \ast r_{\delta }\left( x,t\right)
=\int_{0}^{t}r_{0}\left( t-t^{\prime }\right) r_{\delta }\left( x,t^{\prime
}\right) \,dt^{\prime }.  \label{r(x,t)_convolution}
\end{equation}%
Another relevant initial response is that with Laplace transform:\
\begin{equation}
\tilde{r}_{n}\left( x,s\right) =\frac{\exp \left( -s\,n\left( s\right)
\,x/c\right) }{s\,n\left( s\right) },  \label{r_n_laplace}
\end{equation}%
so $r_{\delta }\left( x,t\right) $ can be obtained from $r_{n}\left(
x,t\right) $ by partial derivative with respect to $x$. Of course, the
corresponding Laplace transforms exhibit the same singular points, i.e.
those of $n\left( s\right) $.

We have two ways to represent the solution (\ref{r(x,t)_inv_Laplace}), that
is versus $x$ (at fixed time $t$), and versus $t$ (at fixed position $x$).
In the first case, we have $0\leq x\leq ct$, while in the second case, we
have $t\geq x/c\geq 0$. Correspondingly we introduce the parameters:%
\begin{equation}
\mu =\frac{x}{ct}\quad \left( 0\leq \mu \leq 1\right) ,\quad \theta =\frac{ct%
}{x}\quad \left( 1\leq \theta \leq \infty \right) .  \label{parameters}
\end{equation}

In order to decrease the computational difficulties for the Laplace inversion,
we propose to deform the original path of integration (the Bromwich path) in
(\ref{r(x,t)_inv_Laplace}) into another equivalent to it (unless possible
contributions of singularities) that is expected to be more convenient, i.e.
the steepest descent path through the saddle points of the complex function:%
\begin{equation}
F_{\mu }\left( s\right) =s\left[ 1-\mu \,n\left( s\right) \right] ,\quad
F_{\theta }\left( s\right) =s\left[ 1-\,n\left( s\right) /\theta \right] ,
\label{F_mu_def}
\end{equation}%
according to our choice of representation in (\ref{parameters}).
We refer to the reader to
\cite{Brillouin,MainardiSteepest}
for a detailed explanation of the steepest descent method and its applications.

\section{The Klein-Gordon equation  with \\ dissipation}

A model equation for uniaxial waves in dispersive media with dissipation is
the Klein-Gordon equation with an additional term that takes into account of
attenuation due to dissipation. We refer to this equation as the
Klein-Gordon with dissipation, i.e. KGD equation \cite{MainardiKlein}:%
\begin{equation}
r_{tt}+a\,r_{t}+b\,r-c^{2}r_{xx}=0,  \label{KGD_def}
\end{equation}%
where $r=r\left( x,t\right) $ is the response variable, $c^{2}$ denotes the
square of the wave-front velocity, and $a,\,b$ are non negative constants.
If $b=0$ the equation reduces to the so-called {telegraph equation},
whereas if $a=0$ we recover the classical  {Klein-Gordon equation}
without dissipation. The space-time coordinates are taken in the quadrant
 $x,\,t\geq 0$. Also, we keep the usual boundary and initial conditions:%
\begin{eqnarray}
r\left( 0,t\right)  &=&r_{0}\left( t\right) ,\quad \lim_{x\rightarrow \infty
}r\left( x,t\right) =0,  \label{Boundary_conds} \\
r\left( x,0\right)  &=&r_{t}\left( x,0\right) =0.  \label{Initial_conds}
\end{eqnarray}%
For the complex index of refraction associated to the KGD  equation
we have\footnote{In the published version the factor $s$ is missed in front of $n(s)$ in the LHS of Eq. (2.4).}  %
\begin{equation}
s\, n\left( s\right) =\frac{\sqrt{s^{2}+as+b}}{s}=\frac{\sqrt{\left(
s+a/2\right) ^{2}+\Delta }}{s},\quad \Delta =b^{2}-\frac{a^2}{4}.  \label{n(s)_KGD}
\end{equation}%
As a consequence, we write the following solutions in the Laplace domain:\
\begin{eqnarray}
\tilde{r}_{\delta }\left( x,s\right)  &=&\exp \left( -\frac{x}{c}\sqrt{%
\left( s+a/2\right) ^{2}+\Delta }\right) ,  \label{r_delta_Laplace_KGD} \\
\tilde{r}_{n}\left( x,s\right)  &=&\frac{\exp \left( -\frac{x}{c}\sqrt{%
\left( s+a/2\right) ^{2}+\Delta }\right) }{\sqrt{\left( s+a/2\right)
^{2}+\Delta }}.  \label{r_n_Laplace_KGD}
\end{eqnarray}%
For the impulse response, we get 
(see \cite[Eqns. (126)-(127), p. 213]{Ghizzetti}):

\begin{itemize}
\item If $\Delta >0$:%
\begin{eqnarray}
r_{\delta }\left( x,t\right) &=&\exp \left( -\frac{a\,x}{2c}\right) \delta
\left( t-\frac{x}{c}\right)  \label{r_delta_DELTA>0} \\
&&-\frac{x}{c}\sqrt{\Delta }\exp \left( -\frac{at}{2}\right) \frac{%
J_{1}\left( \sqrt{\Delta \left( t^{2}-\frac{x^{2}}{c^{2}}\right) }\right) }{%
\sqrt{t^{2}-\frac{x^{2}}{c^{2}}}}\Theta \left( t-\frac{x}{c}\right) ,
\nonumber
\end{eqnarray}
\end{itemize}
\begin{itemize}
\item If $\Delta <0$:%
\begin{eqnarray}
r_{\delta }\left( x,t\right) &=&\exp \left( -\frac{a\,x}{2c}\right) \delta
\left( t-\frac{x}{c}\right)  \label{r_delta_DELTA<0} \\
&&+\frac{x}{c}\sqrt{-\Delta }\exp \left( -\frac{at}{2}\right) \frac{%
I_{1}\left( \sqrt{-\Delta \left( t^{2}-\frac{x^{2}}{c^{2}}\right) }\right) }{%
\sqrt{t^{2}-\frac{x^{2}}{c^{2}}}}\Theta \left( t-\frac{x}{c}\right) ,
\nonumber
\end{eqnarray}
\end{itemize}
where $\Theta \left( t\right) $ denotes the Heaviside theta function, and $%
J_{1}\left( t\right) $ and $I_{1}\left( t\right) $ are the Bessel and
modified Bessel function of first order, respectively \cite[Chaps. 49\&52]{Atlas}.

For the other response, we get 
(see \cite[Eqns. (124)-(125) p. 213]{Ghizzetti}):

\begin{itemize}
\item If $\Delta >0$:%
\begin{equation}
r_{n}\left( x,t\right) =\exp \left( -\frac{at}{2}\right) \frac{J_{0}\left(
\sqrt{\Delta \left( t^{2}-\frac{x^{2}}{c^{2}}\right) }\right) }{\sqrt{t^{2}-%
\frac{x^{2}}{c^{2}}}}\Theta \left( t-\frac{x}{c}\right) ,
\label{r_n_DELTA>0}
\end{equation}
\end{itemize}
\begin{itemize}
\item If $\Delta <0$:%
\begin{equation}
r_{n}\left( x,t\right) =\exp \left( -\frac{at}{2}\right) \frac{I_{0}\left(
\sqrt{-\Delta \left( t^{2}-\frac{x^{2}}{c^{2}}\right) }\right) }{\sqrt{t^{2}-%
\frac{x^{2}}{c^{2}}}}\Theta \left( t-\frac{x}{c}\right) ,
\label{r_n_DELTA<0}
\end{equation}
\end{itemize}
where now $J_{0}\left( t\right) $ and $I_{0}\left( t\right) $ are the Bessel
and modified Bessel function of zeroth order, respectively. In particular
case of the telegraph equation, i.e. $b=0$, we get $\Delta =-a^{4}/4<0$, and
the corresponding response was found by Lee-Kanter for the Maxwell model of
linear viscoelasticity \cite{LeeKanter}.

We wish to approximate the solutions $r_{\delta }\left( x,t\right) $ and $%
r_{n}\left( x,t\right) $, computing the inverse Laplace transforms given in (%
\ref{r_delta_Laplace_KGD})\ and (\ref{r_n_Laplace_KGD}), by means of the
steepest descent method, and then compare the results to the exacts solution
provided in (\ref{r_delta_DELTA<0})-(\ref{r_n_DELTA>0}). This method
consists in replacing the Bromwich path with the so-called \textit{steepest
descent path} (SPD), that requires finding the saddle points of the $F_{\mu
} $ function defined in (\ref{F_mu_def}). For the particular case of the
KGD\ equation, insert (\ref{n(s)_KGD})\ into (\ref{F_mu_def}) to obtain:\
\begin{equation}
F_{\mu }\left( s\right) =s-\mu \,\sqrt{\left( s+a/2\right) ^{2}+\Delta }.
\label{F_mu_KGD}
\end{equation}

We recall that the required steepest descent path (see \cite{Brillouin,MainardiSteepest}) is the path through saddle
points of $F_{\mu}$
 along which the real part of $F_{\mu}$
 attains its maximum value so that the imaginary part of
  $F_{\mu}$
 is constant. We note that the steepest descent path, denoted in the
following as $\gamma_{\mu}$, depends on the parameter $\mu$.

In order to proceed with the determination of the saddle points of 
$F_{\mu}$, which are the solution of $F_{\mu }^{\prime }\left( s\right) =0$, and
the steepest descent path $\gamma_{\mu}$, it is convenient to discuss
separately the cases with $\Delta <0$ and $\Delta >0$. Indeed, we deal at
first the easiest case $\Delta <0$ because the corresponding $\gamma_{\mu}$
turns out to be a closed curve, while for $\Delta >0$ 
the corresponding $\gamma_{\mu}$ turns out to be a couple of open curves symmetric with respect
to the negative real axis.

In the figures of the SDP, the saddle points are characterized by arrows
that show the direction of ascent along the lines of steepest descent,
according the convection used by Brilluoin.

\section{The determination of the steepest descent path}

\subsection{Case $\Delta <0$}

\begin{figure}[h]
	\begin{center}
		\includegraphics[width=0.8\textwidth]{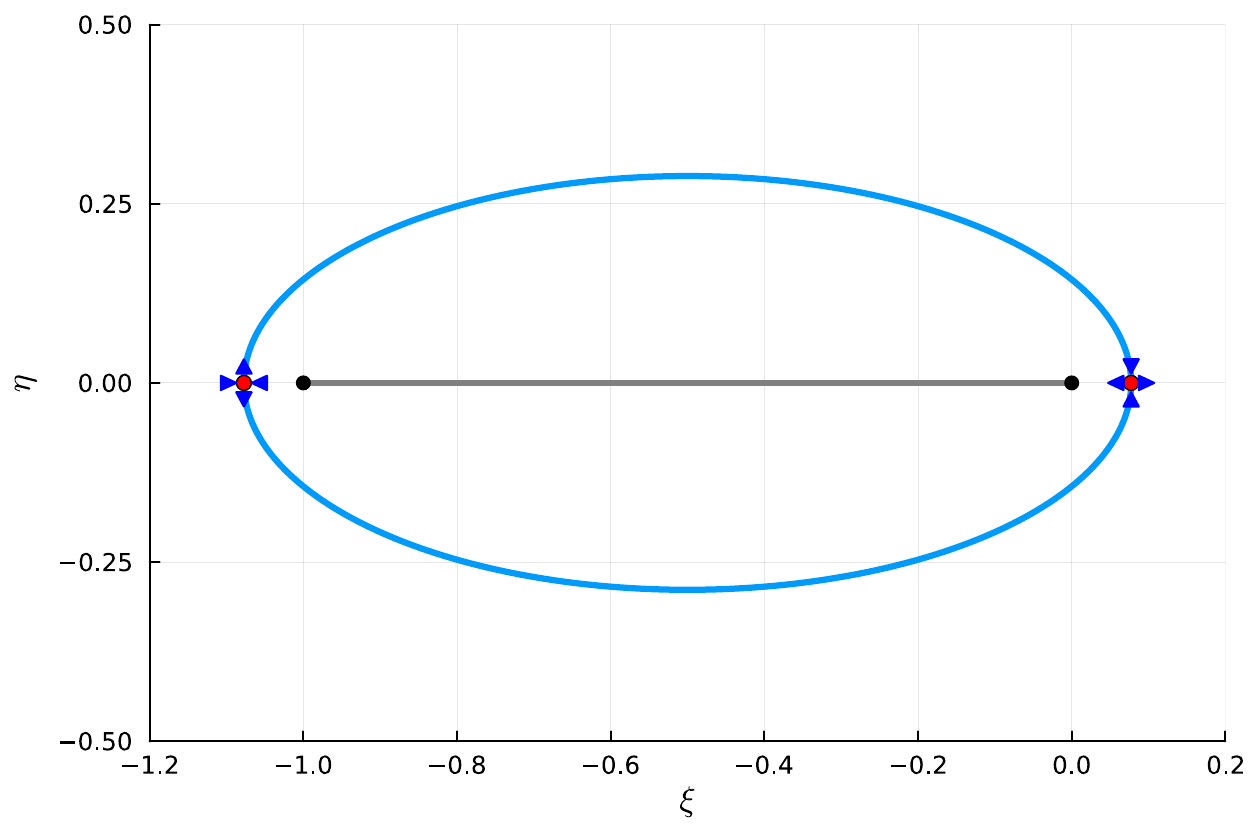}	
	\end{center}
	\vskip -0.5truecm
	\caption{Steepest descent path for the case $\Delta < 0$ (blue curve). The red and black dots represent, respectively, the saddle points and the branch points; the grey line represent the branch cut (values of the parameters: $a=1$, $b=0$, $c=1$, $\Delta=-1/4$; $\mu=1/2$). \label{fig:sdp-delta-neg} }
\end{figure}

If $\Delta <0$, the saddle points of the function $F_{\mu }\left( s\right) $
are the following:%
\begin{equation}
p_{1}=-\frac{a}{2}-\sqrt{\frac{-\Delta }{1-\mu ^{2}}},\quad p_{2}=-\frac{a}{2%
}+\sqrt{\frac{-\Delta }{1-\mu ^{2}}},  \label{p1_p2_DELTA<0}
\end{equation}%
which are both real. The values of $F_{\mu }$ at the saddle points $p_{1}$
and $p_{2}$ are respectively $\phi _{1}=\varphi _{1}+\omega _{1}\,i$ and $%
\phi _{2}=\varphi _{2}+\omega _{2}\,i$ with%
\begin{eqnarray}
\varphi _{1} &=&-\frac{a}{2}+\left( 1-\mu ^{2}\right) \sqrt{\frac{-\Delta }{%
1-\mu ^{2}}},  \label{phi_1_DELTA<0} \\
\varphi _{2} &=&-\frac{a}{2}-\left( 1+\mu ^{2}\right) \sqrt{\frac{-\Delta }{%
1-\mu ^{2}}},  \label{phi_2_DELTA<0} \\
\omega _{1} &=&\omega _{2}=0.  \label{w1=w2_0_DELTA<0}
\end{eqnarray}%
Since $\phi _{1}$ and $\phi _{2}$ are both reals (i.e. the imaginary part of
$F_{\mu }$ is zero in the two saddle points) the steepest descent path is in
this case a curve in the complex plane passing through the two saddle points
$p_{1}$ and $p_{2}$. In passing, we note that the branch points (namely, the
points for which the argument of the square root vanishes) are the following:%
\begin{equation}
b_{1}=-\frac{a}{2}-\sqrt{-\Delta },\quad b_{2}=-\frac{a}{2}+\sqrt{-\Delta }.
\label{b1_b2_DELTA<0}
\end{equation}%
The branch points, as the saddle points, are on the real axis, and
\begin{equation}
p_{1}\leq b_{1}<-\frac{a}{2}<b_{2}\leq p_{2},  \label{Inequation_DELTA<0}
\end{equation}%
where the equal signs hold for the limit value $\mu=0$.
\newpage
On the steepest descent path though the saddle points $p_{1}$, $p_{2}$, the
imaginary part of $F_{\mu }$ is constant (and equal to $\omega _{1}=\omega
_{2}=0$), therefore the steepest descent path is found as follows:
\begin{equation}
\mathrm{Im}\left( F_{\mu }\right) =
\mathrm{Im}\left( s-\mu \,\sqrt{\left(
s+a/2\right) ^{2}+\Delta }\right) =0.  \label{Im(F_mu)=0_DELTA<0}
\end{equation}
Setting $s=\xi +\eta \,i$, and $\sqrt{\left( s+a/2\right) ^{2}+\Delta}
=A+B\,i$, with $\xi ,\eta ,A,B\in \mathbb{R}$, we have%
\begin{eqnarray}
\eta -\mu \,B &=&0,  \label{Ellipse_1} \\
\left( \xi +a/2\right) ^{2}+\Delta  &=&A^{2}-B^{2},  \label{Ellipse_2} \\
\left( \xi +\frac{a}{2}\right) \eta  &=&AB,  \label{Ellipse_3}
\end{eqnarray}%
from which we readily find%
\begin{equation}
\frac{\left( \xi +a/2\right) ^{2}}{\alpha ^{2}}+\frac{\eta ^{2}}{\beta ^{2}}%
=1,  \label{Ellipse_DELTA<0}
\end{equation}%
with%
\begin{equation}
\alpha =\sqrt{\frac{-\Delta }{1-\mu ^{2}}},\quad \beta =\mu \alpha =\mu
\sqrt{\frac{-\Delta }{1-\mu ^{2}}}.  \label{semiaxes_ellipse}
\end{equation}%
The steepest descent path $\gamma _{\mu }$ turns out to be a one-parameter
family of ellipses with semi-axes $\alpha $ and $\beta $, 
see Fig. \ref{fig:sdp-delta-neg}
\vskip -0.5truecm

\newpage

\subsection{Case $\Delta >0$}

In this case, the saddle points of the $F_{\mu }$\ function are the
following:%
\begin{equation}
p_{1}=-\frac{a}{2}-i\sqrt{\frac{\Delta }{1-\mu ^{2}}},\quad p_{2}=-\frac{a}{2%
}+i\sqrt{\frac{\Delta }{1-\mu ^{2}}},  \label{p1_p2_DELTA>0}
\end{equation}%
which are complex conjugates. The values of $F_{\mu }$ at the saddle points $%
p_{1}$ and $p_{2}$ are respectively $\phi _{1}=\varphi _{1}+\omega _{1}\,i$
and $\phi _{2}=\varphi _{2}+\omega _{2}\,i$ with%
\begin{eqnarray}
\varphi _{1} &=&\varphi _{2}=-\frac{a}{2},  \label{phi_1_phi_2_DELTA>0} \\
\omega _{1} &=&+\sqrt{\Delta \left( 1-\mu ^{2}\right) },
\label{w1_DELTA>0} \\
\omega _{2} &=&-\sqrt{\Delta \left( 1-\mu ^{2}\right) }.  \label{w2_DELTA>0}
\end{eqnarray}%
The branch points (which are complex conjugates) are the following:%
\begin{equation}
b_{1}=-\frac{a}{2}-i\sqrt{\Delta },\quad b_{2}=-\frac{a}{2}+i\sqrt{\Delta }.
\label{b1_b2_DELTA>0}
\end{equation}%
We note that the saddle points, as well as the branch points, lie on the
line $\xi =-a/2$ parallel to the imaginary axis, and their imaginary parts
are such that%
\begin{equation}
\mathrm{Im}\,p_{1}\leq \mathrm{Im}\,b_{1}<0<\mathrm{Im}\,b_{2}\leq \mathrm{Im%
}\,p_{2},  \label{Inequation_DELTA>0}
\end{equation}%
where the equal signs hold for the limit value
 $\mu=0$.

In contrast to the case with $\Delta<0$, the steepest descent path 
$\gamma_{\mu}$ is now made of two branches: 
a branch through the saddle point
$p_{1}$, denoted 
as $\gamma_{\mu,1}$,
 and a branch through the saddle point $p_{2}$, denoted as 
 $\gamma _{\mu,2}$.

\begin{figure}[htbp]
\begin{center}
\includegraphics[width=0.8\textwidth]{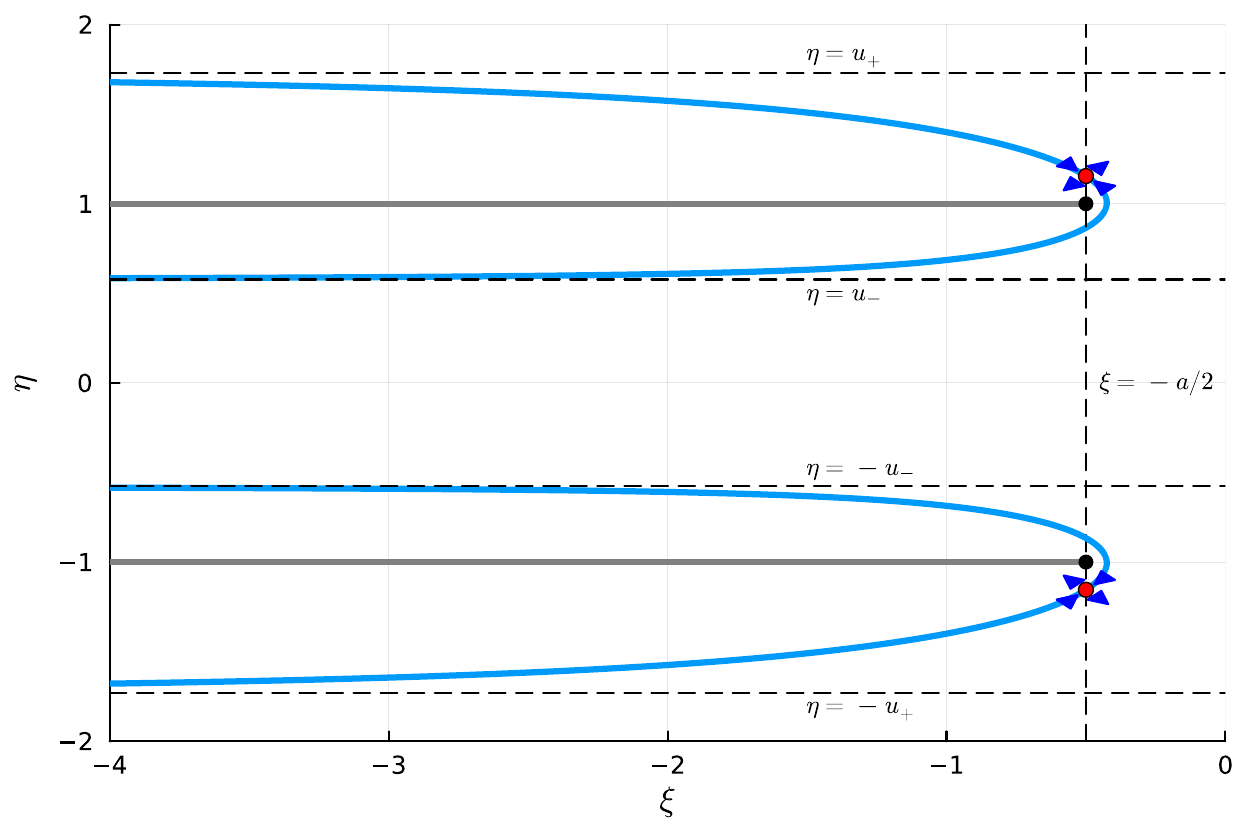}
\end{center}
\vskip -0.5truecm
\caption{Steepest descent path for the case $\Delta > 0$ (blue curve). The
red and black dots represent, respectively, the saddle points and the branch
points; the grey line represent the branch cut (values of the parameters: $%
a=1$, $b=5/4$, $c=1$, $\Delta=1$; $\protect\mu=1/2$). }
\label{fig:sdp-delta-pos}
\end{figure}

Along the branch of the steepest descent path though the saddle point $p_{k}$
($k=1,2$), the imaginary part of $F_{\mu}$ is constant and equal to 
$\omega _{k}$, i.e.%
\begin{equation}
\mathrm{Im}\left( F_{\mu }\right) =\mathrm{Im}\left( s-\mu \,\sqrt{\left(
s+a/2\right) ^{2}+\Delta }\right) =\omega _{k}.
\label{Im_(F_mu)=w_k_DELTA>0}
\end{equation}%
Proceeding as in the case with $\Delta <0$, we find%
\begin{eqnarray}
\eta -\mu \,B &=&\omega _{k},  \label{g_k_curve_1} \\
\left( \xi +a/2\right) ^{2}-\eta ^{2}+\Delta  &=&A^{2}-B^{2},
\label{g_k_curve_2} \\
\left( \xi +\frac{a}{2}\right) \eta  &=&AB,  \label{g_k_curve_3}
\end{eqnarray}%
from which we find for $k=1,2$
\begin{equation}
\left( \xi +\frac{a}{2}\right) ^{2}=\left( \frac{\eta -\omega _{k}}{\mu }%
\right) ^{2}\frac{\left( \eta -\omega _{k}\right) ^{2}+\mu ^{2}\left( \Delta
-\eta ^{2}\right) }{\mu ^{2}\eta ^{2}-\left( \eta -\omega _{k}\right) ^{2}},
\label{gamma_DELTA>0}
\end{equation}%
which defines the branch 
$\gamma _{\mu,1}$ of the steepest descent path through the saddle point 
$p_{1}$ (for $k=1$), and the branch
 $\gamma _{\mu,2}$ through the saddle point $p_{2}$ (for $k=2$). 
 The steepest descent path
is defined as
$\gamma_{\mu,1}\cup \gamma _{\mu,2}$. 
An analysis of (\ref{gamma_DELTA>0}) reveals that the two branches of
the path $\gamma _{\mu}$, i.e. 
$\gamma _{\mu,1}$ and $\gamma _{\mu,2}$
 are symmetric with respect to the real axis, see
  Fig \ref{fig:sdp-delta-pos}.
%%%
%%%
An inspection of the properties of $\gamma _{\mu,k}$
shows that each branch of $\gamma _{\mu,k}$
 has two horizontal asymptotes given by:
\[
\eta =u_{\pm }\equiv \sqrt{\frac{1\pm \mu }{1\mp \mu }\,\Delta },
\]
where, taking as parameter $u=\eta $, we have the parametrization:
\begin{equation}
\gamma _{\mu,k}
\left( u\right) =g_{k}\left( u\right) +i\,u,\quad u_{-}<(-1)^{k}u<u_{+},
\label{gamma_mu_k_def}
\end{equation}
where
\begin{equation}
g_{k}\left( u\right) =-\frac{a}{2}\pm \frac{\left\vert u-\omega
_{k}\right\vert }{\mu }\sqrt{\frac{\left( u-\omega _{k}\right) ^{2}+\mu
^{2}\left( \Delta -u^{2}\right) }{\mu ^{2}u^{2}-\left( u-\omega _{k}\right)
^{2}}}.  \label{g_k_def}
\end{equation}

\section{Numerical evaluation of $r_{\protect\delta }\left( x,t\right) $}

\subsection{Case $\Delta <0$}

\begin{figure}[htbp]
\begin{center}
\includegraphics[width=0.8\textwidth]{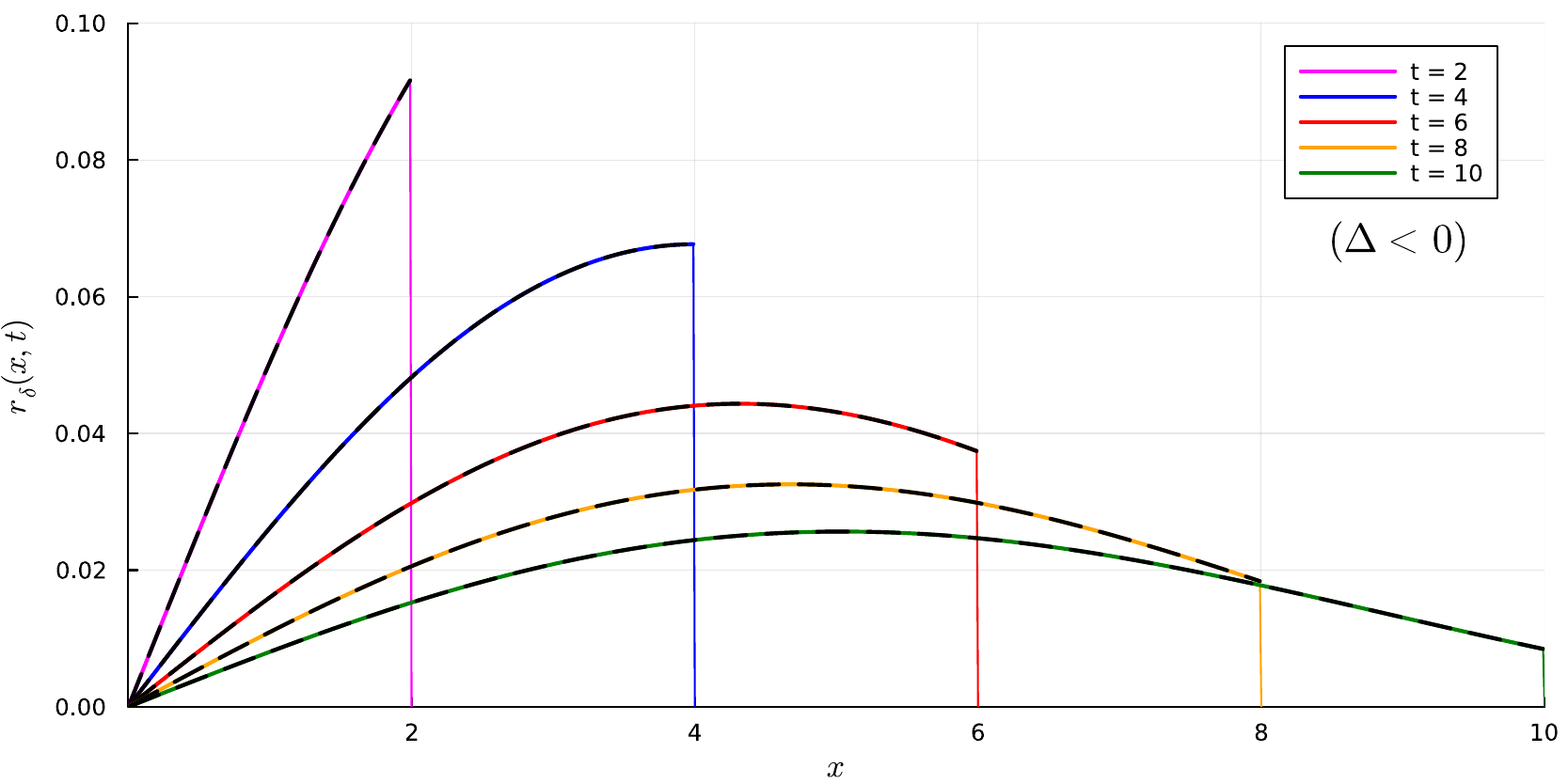}
\end{center}
\vskip -0.5truecm
\caption{Comparison between exact and numerical solutions
of $r_\delta(x,t)$ for  %%\% the case 
$\Delta < 0$, for several values of $t$ (values of the parameters: $a=1$, $b=0
$, $c=1$, $\Delta=-1/4$). }
\label{fig:sol-delta-neg}
\end{figure}

In this case, $\gamma _{\mu }$ is a closed curve which can be easily
parameterized as follows:\
\begin{equation}
\gamma _{\mu }\left( u\right) =-\frac{a}{2}+\alpha \cos u+i\,\beta \sin
u,\qquad 0\leq u<2\pi ,  \label{gamma_mu_DELTA<0}
\end{equation}%
from which:%
\begin{equation}
\gamma _{\mu }^{\prime }\left( u\right) =-\alpha \sin u+i\,\beta \cos u.
\label{gamma'(u)_def}
\end{equation}%
Letting%
\begin{equation}
f_{\delta }\left( x,t;s\right) =\frac{1}{2\pi \,i}\exp \left( st-\frac{x}{c}%
\sqrt{\left( s+\frac{a}{2}\right) ^{2}+\Delta }\right) ,  \label{f_delta_def}
\end{equation}%
with the above parametrization, the integral along $\gamma _{\mu }$ turns
into a real line integral (of the complex $f_{\delta }$ function) as
follows:\
\begin{equation}
r_{\delta }\left( x,t\right) =\int_{\gamma _{\mu }}f_{\delta }\left(
x,t;s\right) \,ds=\int_{0}^{2\pi }f_{\delta }\left( x,t;\gamma _{\mu }\left(
u\right) \right) \,\gamma _{\mu }^{\prime }\left( u\right) \,du.
\label{r_delta_steepest_DELTA<0}
\end{equation}%
The resulting integral can be numerically evaluated by means of standard
numerical method, such as the adaptative Gauss-Kronrod quadrature.
 A comparative between numerical and exact results is shown in 
 Fig. \ref{fig:sol-delta-neg}, where the solution $r_{\delta }\left( x,t\right) $ is
plotted as a function of $x$ for several values of $t$.
%%%%
%%%%
\newpage
\subsection{Case $\Delta >0$}
\begin{figure}[htbp]
\begin{center}
\includegraphics[width=0.8\textwidth]{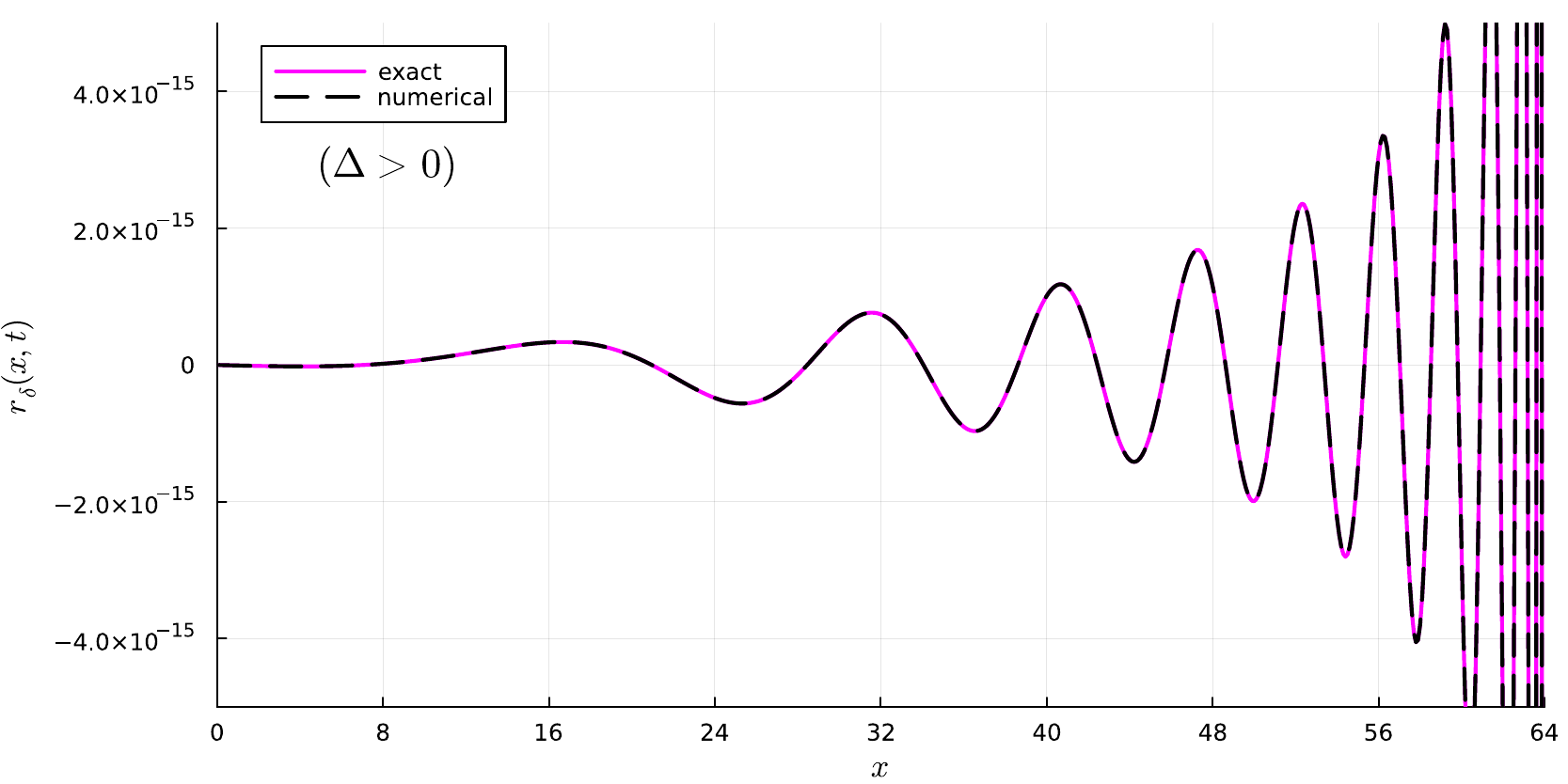}
\end{center}
\vskip -0.5truecm
\caption{Comparison between exact and numerical results 
of $r_\delta(x,t)$ for % the case with 
$\Delta > 0$ (values of the parameters: $a=1$, $b=5/4$, $c=1$, $\Delta=1$, $%
t=64$). }
\label{fig:sol-delta-pos-1}
\end{figure}

\begin{figure}[htbp]
\begin{center}
\includegraphics[width=0.8\textwidth]{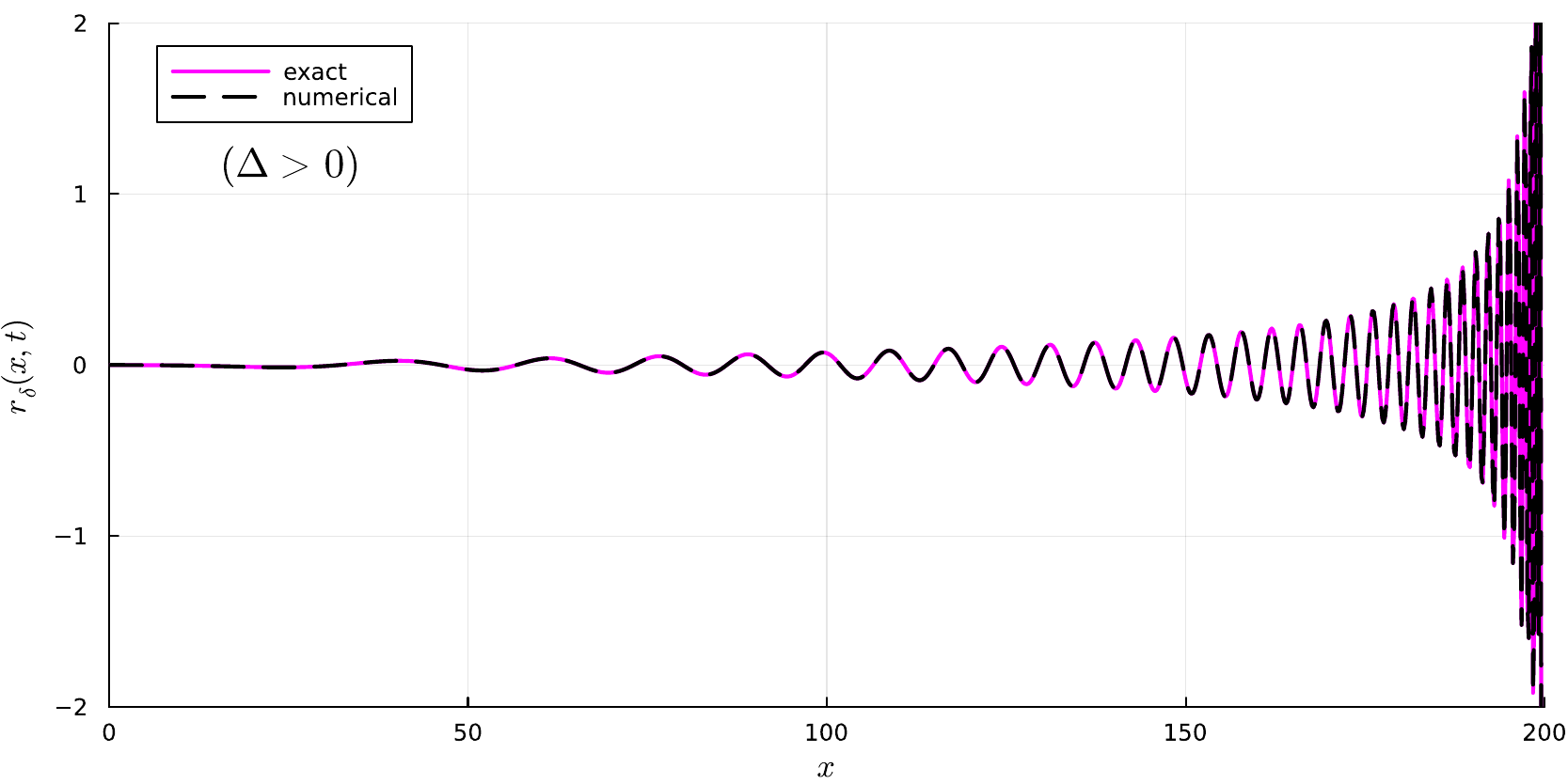}
\end{center}
\vskip -0.5truecm
\caption{Comparison between exact and numerical results
of $r_\delta(x,t)$  for % the case with 
$ \Delta > 0$ (values of the parameters: $a=10^{-4}$, $b=5$, $c=2$, $%
\Delta=4.9999999975$, $t=100$). }
\label{fig:sol-delta-pos-2}
\end{figure}

%%%
Following the procedure outlined for the case $\Delta <0$, and replacing the
Bromwich path by the descendent $\gamma _{\mu }$, we have%
\begin{equation}
r_{\delta }\left( x,t\right) =\int_{\gamma _{\mu ,1}\cup \gamma _{\mu
,2}}f_{\delta }\left( x,t;s\right) \,ds=\mathcal{I}_{\delta ,1}+\mathcal{I}%
_{\delta ,2},  \label{r_delta_steepest_DELTA>0}
\end{equation}%
where $\forall k=1,2$
\begin{equation}
\mathcal{I}_{\delta ,k}=\int_{\gamma _{\mu ,k}}f_{\delta }\left(
x,t;s\right) \,ds.  \label{I_delta,k_def}
\end{equation}%
Deploying the parametrization discussed above, we have%
\begin{equation}
\mathcal{I}_{\delta ,k}=\left( -1\right) ^{k}\int_{\left( -1\right)
^{k}u_{-}}^{\left( -1\right) ^{k}u_{+}}\,f_{\delta }\left( x,t;\,\gamma _{%
%TCIMACRO{\U{3bc} }%
%BeginExpansion
\mu
%EndExpansion
,k}\left( u\right) \right) \,\gamma _{%
%TCIMACRO{\U{3bc} }%
%BeginExpansion
\mu
%EndExpansion
,k}^{\prime }\left( u\right) \,du,  \label{I_delta_k_def}
\end{equation}%
where $\gamma _{%
%TCIMACRO{\U{3bc} }%
%BeginExpansion
\mu
%EndExpansion
,k}$ is given in (\ref{gamma_mu_k_def})-(\ref{g_k_def}), and $\gamma _{%
%TCIMACRO{\U{3bc} }%
%BeginExpansion
\mu
%EndExpansion
,k}^{\prime }\left( u\right) =g_{k}^{\prime }\left( u\right) +i$. As for the
case with $\Delta <0$, the Bromwich integral is therefore replaced by real
line integrals, which can be evaluated by means of standard algorithms. In
Figs. \ref{fig:sol-delta-pos-1}\ and \ref{fig:sol-delta-pos-2}, comparisons
between numerical and exact results are plotted for two set of parameter
values.

\subsection{Useful property of the integrals for $r_{\protect\delta }\left(
x,t\right) $}

Since the exact solutions given in (\ref{r_delta_DELTA>0}) and (\ref%
{r_delta_DELTA<0}) are clearly real, and the solutions given in (\ref%
{r_delta_steepest_DELTA>0}) and (\ref{r_delta_steepest_DELTA<0}) are
respectively equivalent to them, then (\ref{r_delta_steepest_DELTA>0}) and (%
\ref{r_delta_steepest_DELTA<0}) must be real as well. Therefore, we expect
that the real line integrals (of the complex function $f_{\delta }$)
emerging from the application of the steepest descent method must be real.

On this regard, we observe that the steepest descent paths are, in both
cases with $\Delta <0$ and $\Delta >0$, symmetric with respect to the real
axis, see Figs. \ref{fig:sdp-delta-neg}\ and \ref{fig:sdp-delta-pos}. For
instance, for the case $\Delta <0$, we can split the integral given in (\ref%
{r_delta_steepest_DELTA<0})\ as
\begin{eqnarray}
&&r_{\delta }\left( x,t\right)   \label{r_delta_split} \\
&=&\underset{\mathcal{H}_{1,\delta }}{\underbrace{\int_{0}^{\pi }f_{\delta
}\left( x,t;\gamma _{\mu }\left( u\right) \right) \,\gamma _{\mu }^{\prime
}\left( u\right) \,du}}+\underset{\mathcal{H}_{2,\delta }}{\underbrace{%
\int_{\pi }^{2\pi }f_{\delta }\left( x,t;\gamma _{\mu }\left( u\right)
\right) \,\gamma _{\mu }^{\prime }\left( u\right) \,du}}.  \nonumber
\end{eqnarray}%
We appreciate that, due to the symmetry properties of the paths and the
features of the integrand function $f$, the integrals $\mathcal{H}_{1,\delta
}$ and $\mathcal{H}_{2,\delta }$ are complex conjugates, i.e. $\mathcal{H}%
_{2,\delta }=\mathcal{\bar{H}}_{1,\delta }$, which allows us to write%
\begin{equation}
r_{\delta }\left( x,t\right) =2\,\mathrm{Re}\left( \mathcal{H}_{1,\delta
}\right) .  \label{r_delta=2*real}
\end{equation}%
This observation is beneficial because it guarantees to obtain a real
solution without any spurious imaginary part that could emerge from the
numerical approximation of the integrals $\mathcal{H}_{1}$
 and $\mathcal{H}_{2}$.

\begin{figure}[htbp]
\begin{center}
\includegraphics[width=0.7\textwidth]{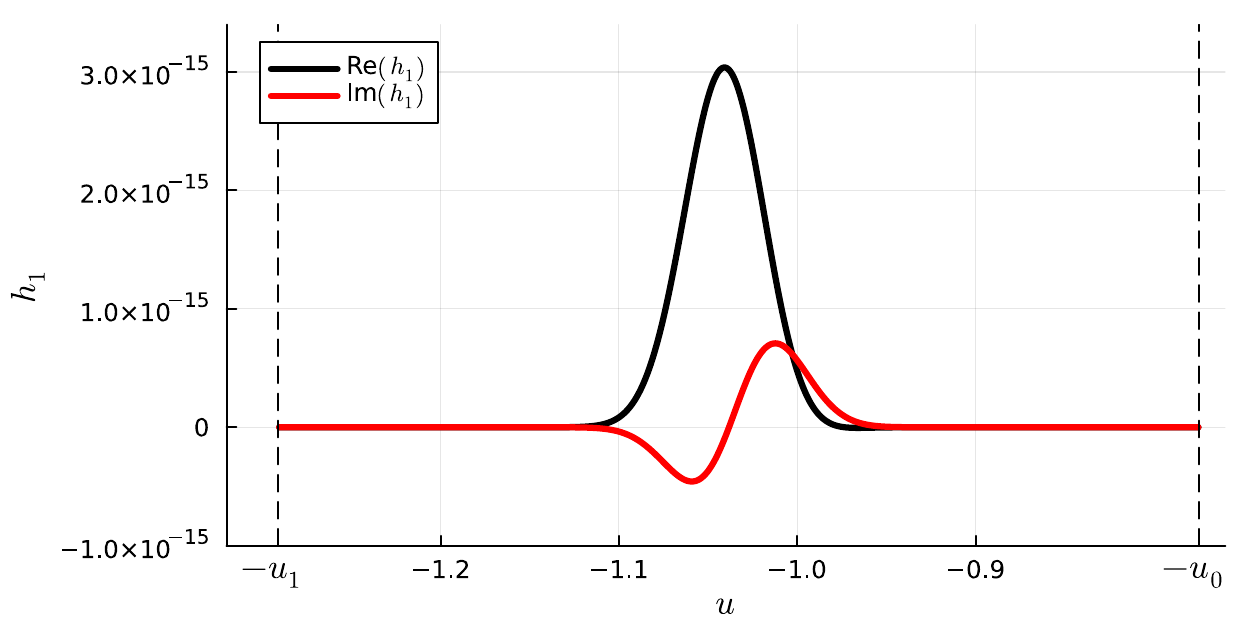} %
\includegraphics[width=0.7\textwidth]{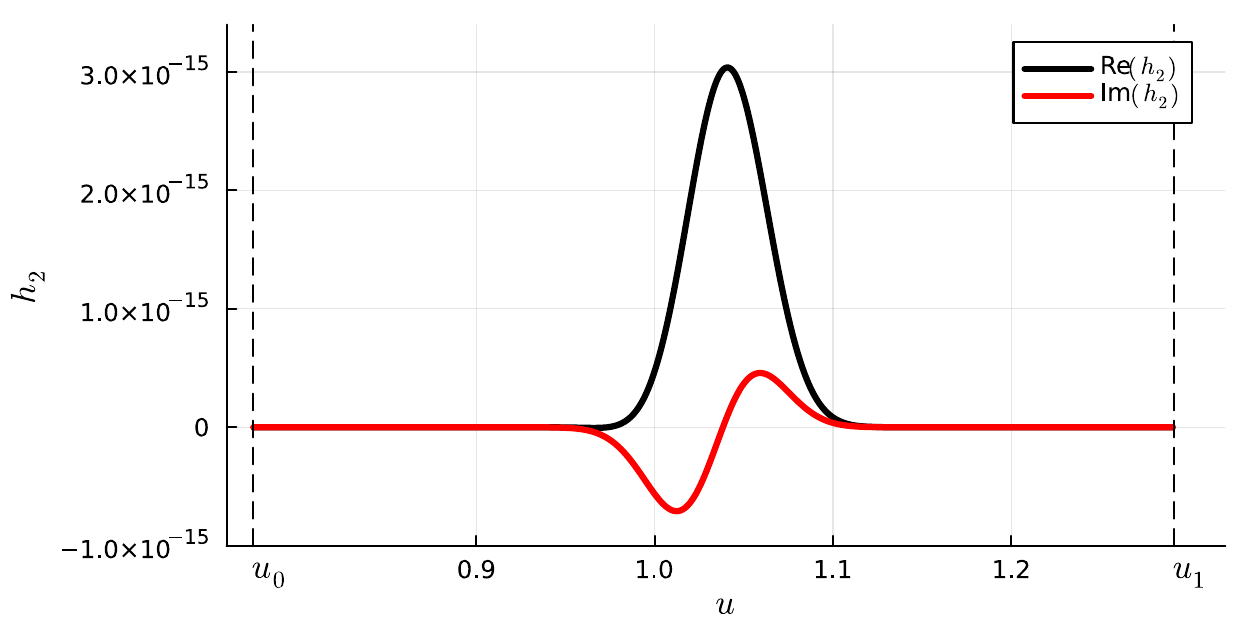}
\end{center}
\caption{Real and imaginary parts of the functions $h_k\left(x, t; u\right)
= f_{\protect\delta}\left(x, t; \protect\gamma_{\protect\mu%
,k}\left(u\right)\right) \protect\gamma_{\protect\mu,k}^{\prime}\left(u%
\right)$ ($k = 1, 2$) as functions of the variable of integration $u$
(values of the parameters: $a=1$, $b=5/4$, $c=1$, $\Delta=1$, $t=64$, $x=16$%
). }
\label{fig:sol-delta-pos-h}
\end{figure}

In Fig. \ref{fig:sol-delta-pos-h}, the integrand functions 
$h_{k}\left(x,t;u\right) =f_{\delta }\left( x,t;\, \gamma _{\mu,k}
\left( u\right) \right) \, \gamma _{\mu,k}^{\prime }\left( u\right) $
 appearing in (\ref{I_delta_k_def}) for an exemplary case are plotted.

\newpage
\section{Numerical evaluation of $r_{n}\left( x,t\right) $}

\subsection{Case $\Delta <0$}
\begin{figure}[htbp]
\begin{center}
\includegraphics[width=0.8\textwidth]{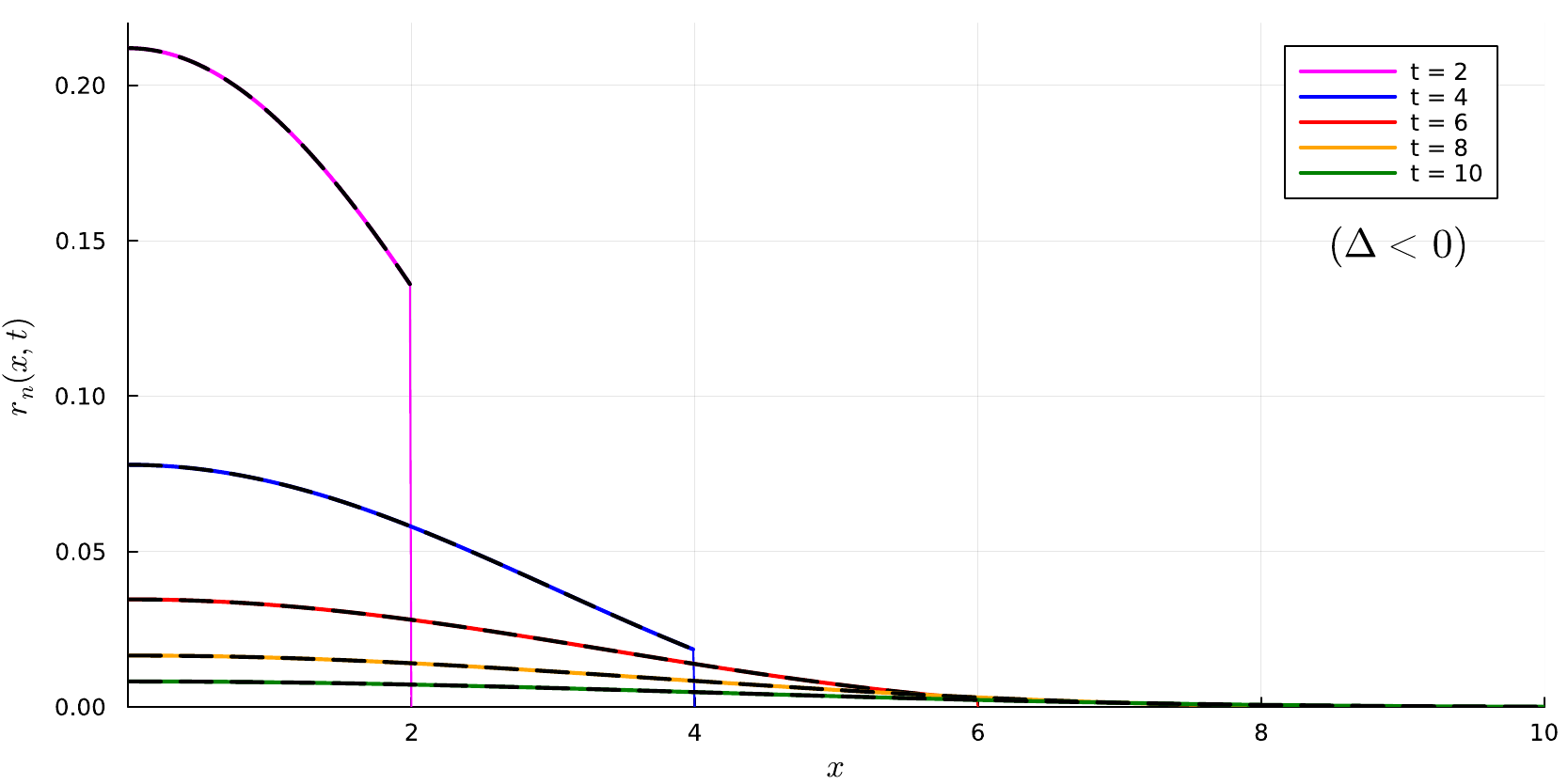}
\end{center}
\vskip -0.5truecm
\caption{Comparison between exact and numerical solutions
of $r_n(x,t)$  for the case $%
\Delta < 0$, for several values of $t$ (values of the parameters: $a=2$, $%
b=1/2$, $c=1$, $\Delta=-1/2$). }
\label{fig:sol-n-neg}
\end{figure}

%%%
In this case, the parametrization $\gamma _{\mu }$ is given in (\ref%
{gamma_mu_DELTA<0}). Letting%
\begin{equation}
f_{n}\left( x,t;s\right) =\frac{1}{2\pi \,i}\frac{\exp \left( st-\frac{x}{c}%
\sqrt{\left( s+a/2\right) ^{2}+\Delta }\right) }{\sqrt{\left( s+a/2\right)
^{2}+\Delta }},  \label{f_n_def}
\end{equation}%
the integral along $\gamma _{\mu }$ turns out to be a real line integral (of
the complex function $f_{n}$), so that
\begin{equation}
r_{n}\left( x,t\right) =\int_{\gamma _{\mu }}f_{n}\left( x,t;s\right)
\,ds=\int_{0}^{2\pi }f_{n}\left( x,t;\gamma _{\mu }\left( u\right) \right)
\,\gamma _{\mu }^{\prime }\left( u\right) \,du.  
\label{r_n_steepest_DELTA<0}
\end{equation}%
Again, the resulting integral may be numerically evaluated by means of a
standard numerical method, such as the adaptive Gauss-Kronrod quadrature.
 A comparison between numerical and exact results is shown in 
 Fig. \ref{fig:sol-n-neg}, where the solution $r_{n}\left( x,t\right) $ is plotted
as a function of $x$ for several values of $t$.
%%%
%%%%
\subsection{Case $\Delta >0$}
\begin{figure}[htbp]
\begin{center}
\includegraphics[width=0.8\textwidth]{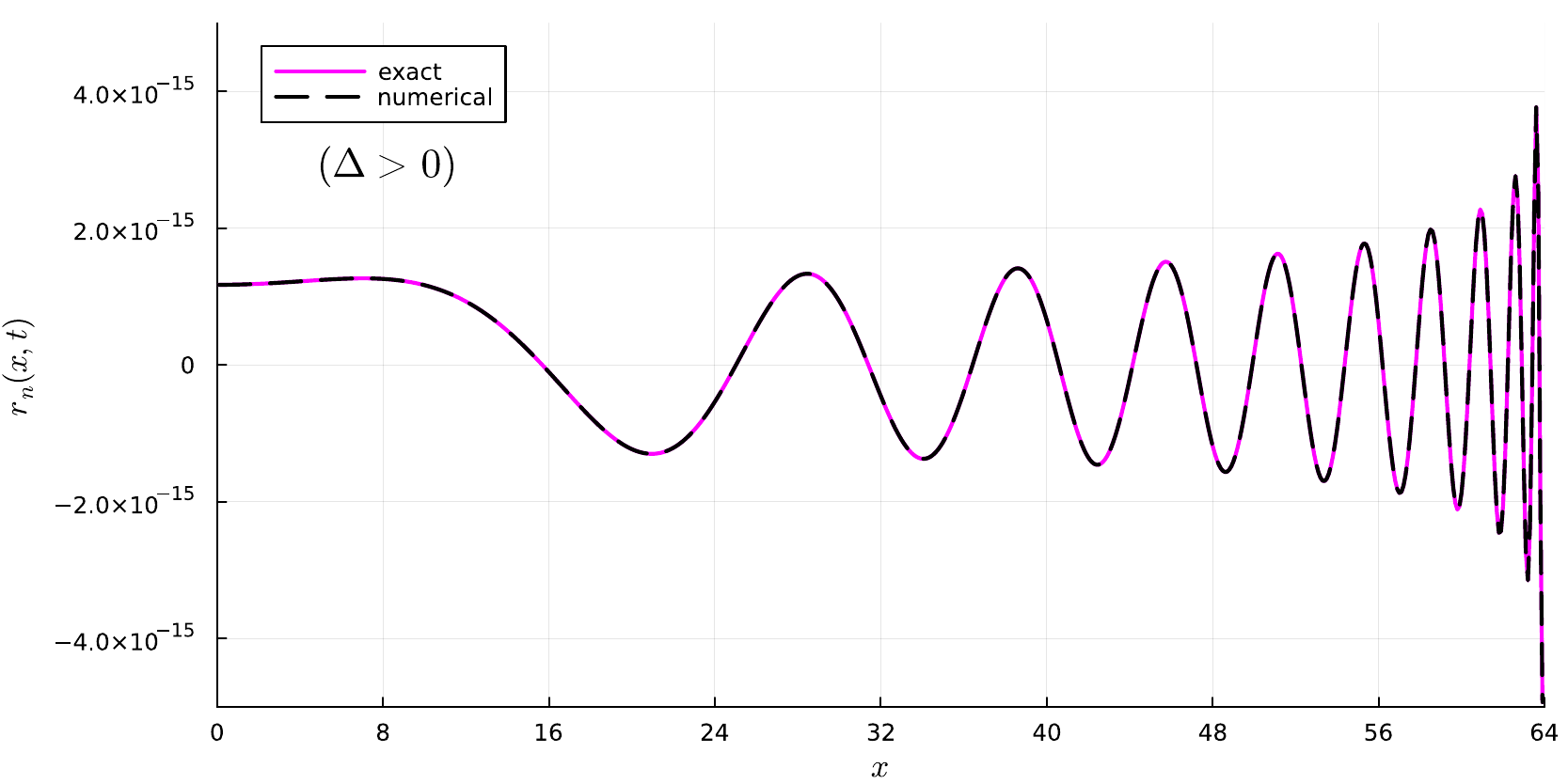}
\end{center}
\vskip -0.5truecm
\caption{ Comparison between exact and numerical results
of  $r_n(x,t)$  for the case with $%
\Delta > 0$ (values of the parameters: $a=1$, $b=5/4$, $c=1$, $\Delta=1$, $%
t=64$). }
\label{fig:sol-n-pos-1}
\end{figure}

\begin{figure}[htbp]
\begin{center}
\includegraphics[width=0.8\textwidth]{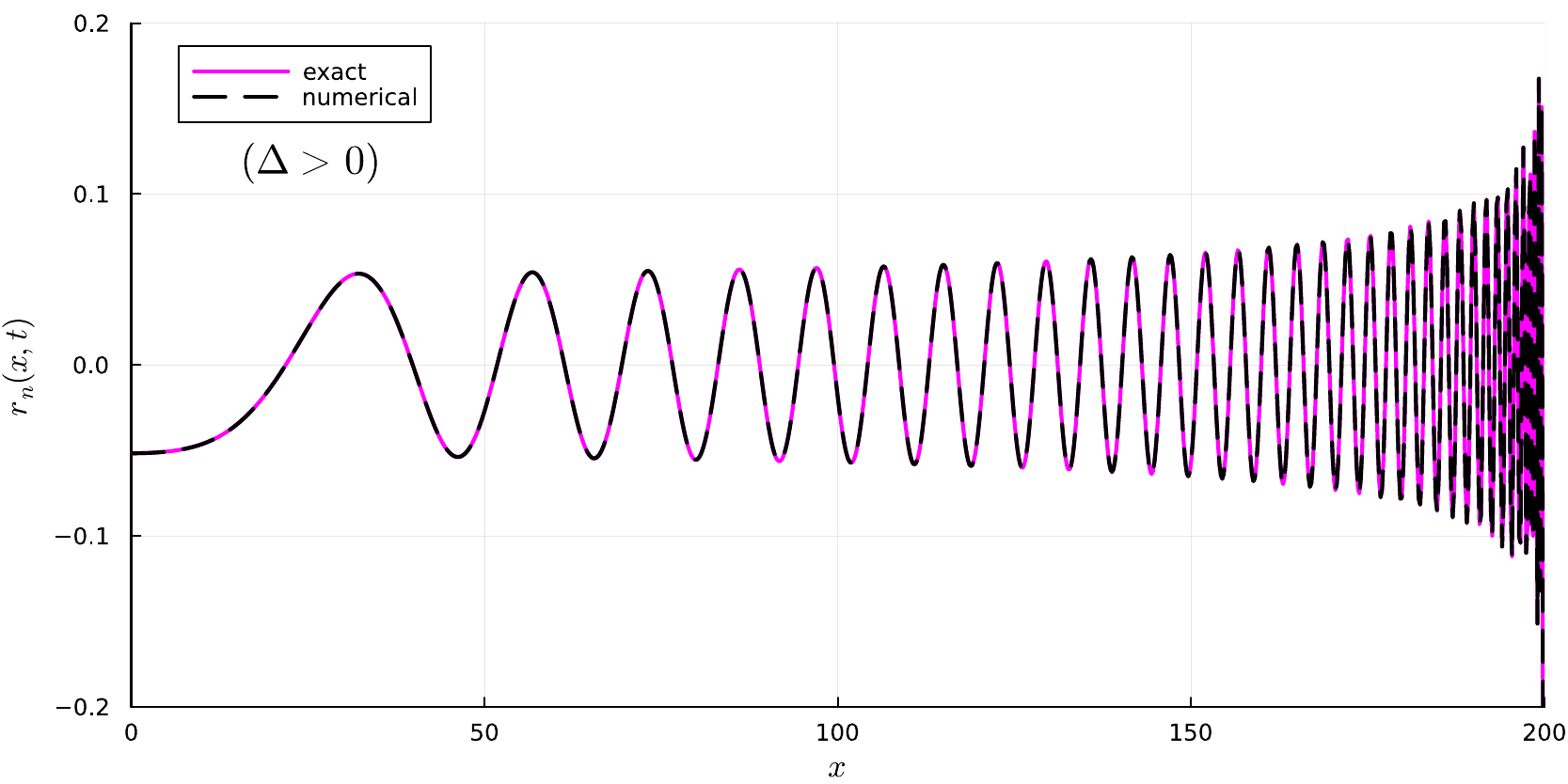}
\end{center}
\vskip -0.5truecm
\caption{Comparison between exact and numerical results
of $r_n(x,t)$  for the case with $%
\Delta > 0$ (values of the parameters: $a=10^{-4}$, $b=5$, $c=2$, $%
\Delta=4.9999999975$, $t=100$). }
\label{fig:sol-n-pos-2}
\end{figure}

%%%
In this case, the parametrization $\gamma _{\mu ,k}$ is given in (\ref%
{gamma_mu_k_def}). Consequently, we have
\begin{equation}
r_{n}\left( x,t\right) =\int_{\gamma _{\mu ,1}\cup \gamma _{\mu
,2}}f_{n}\left( x,t;s\right) \,ds=\mathcal{I}_{n,1}+\mathcal{I}_{n,2},
\label{r_n_steepest_DELTA>0}
\end{equation}%
where $\forall k=1,2$
\begin{equation}
\mathcal{I}_{n,k}=\int_{\gamma _{\mu ,k}}f_{n}\left( x,t;s\right) \,ds.
\label{I_n,k_def}
\end{equation}%
Deploying the parametrization discussed above, we have
\begin{equation}
\mathcal{I}_{n,k}=
\left( -1\right) ^{k}
\int_{\left( -1\right)^{k}u_{-}}^{\left( -1\right) ^{k}u_{+}}\,f_{n}\left( x,t;\,\gamma _{\mu,k}\left( u\right) \right) \,
\gamma _{\mu,k}^{\prime }\left( u\right) \,du,  \label{I_n_k_def}
\end{equation}%
where $\gamma _{\mu,k}$ is given in 
(\ref{gamma_mu_k_def})-(\ref{g_k_def}), 
and 
$\gamma _{\mu,k}^{\prime }\left( u\right) =g_{k}^{\prime }\left( u\right) +i$. 
As for the case with $\Delta <0$, the Bromwich integral is therefore replaced by real line integrals, which can be evaluated by means of standard algorithms. 
In Fig. \ref{fig:sol-n-pos-1} and Fig. \ref{fig:sol-n-pos-2},
comparisons between numerical and exact results are plotted for two set of
parameter values.
%\newpage
%%
%%%%
\newpage
\subsection{Useful property of the integrals for $r_{n}\left( x,t\right) $}

Since the exact solutions given in (\ref{r_n_DELTA>0}) and (\ref{r_n_DELTA<0}) are clearly real, and the solutions given in (\ref{r_n_steepest_DELTA>0})
and (\ref{r_n_steepest_DELTA<0}) are respectively equivalent to them, 
then (\ref{r_n_steepest_DELTA>0}) and (\ref{r_n_steepest_DELTA<0}) must be real as well. 
Therefore, we expect that the real line integrals (of the complex
function $f_{n}$) emerging from the application of the steepest descent
method must be real.

Consequently, analogously to the result obtained previously, it turns out
that
\begin{equation}
r_{n}\left( x,t\right) =2\,\mathrm{Re}\left( \mathcal{H}_{1,n}\right) ,
\label{r_n=2*real}
\end{equation}%
where
\begin{equation}
\mathcal{H}_{1,n}:=\int_{0}^{\pi }f_{n}\left( x,t;\gamma _{\mu }\left(
u\right) \right) \,\gamma _{\mu }^{\prime }\left( u\right) \,du.
\label{H_1,n_def}
\end{equation}

%%%%%%%%%%%%%%%%%%%%%%%%%%%%%%%%%%%%
\newpage
\section{Conclusions}

It is worth stressing that the method discussed here consists in turning the
Bromwich integral into a real line integral of a well-behaved function, for
the solution of which well established (and well performing) quadrature
formulas may be used.

Consequently, unlike the traditional usage of the steepest descent method to calculate the
asymptotic behaviour of the Bromwich integral nearby the saddle points, here
we consider the entire steepest descent path. Therefore, our method does not
provide an asymptotic result, but an exact one, neglecting
the computational error in the numerical evaluation of the corresponding real line integral.

The possibly hard problem of finding $r_{\delta }\left( x,t\right) $ 
and $r_n(x,t)$ is
therefore transferred from the algorithm for the numerical approximation of
the corresponding Bromwich integrals, to the analytical calculation (and parametrization) of the steepest descent path. This is in contrast to the majority of the most popular algorithms available for the inversion of the Laplace
transform, which focus on the numerical solution of the Bromwich integral
(or other equivalent integrals in the complex plane), which is a notoriously
difficult task. Even the most popular algorithms for the numerical inversion
of Laplace transforms (see \cite{Talbot,Weideman2006,Weideman2007}) are known to provide, in critical cases, solutions
affected by very large errors, which --when the exact solution is not
available -- can be difficult to spot, potentially leading to unreliable
results.

\section*{Acknowledgments}

 The research activity of F. Mainardi and A. Mentrelli  has been carried out in the framework of the activities of the National Group of Mathematical Physics
 (GNFM, INdAM). 
\\ 
 The authors are grateful to the anonymous referees for valuable suggestions which help us to improve the presentation of the results.

\section*{Funding}
A. Mentrelli is partially funded by the European Union – NextGenerationEU under the National Recovery and Resilience Plan (PNRR) - Mission 4 Education and research, Component 2 From research to business – Investment 1.1 Notice PRIN 2022 – DD N. 104 dated 2/2/2022, entitled ‘‘The Mathematics and Mechanics of Non-linear Wave Propagation in Solids’’ (proposal code: 2022P5R22 A; CUP: J53D23002350006), and by the Italian National Institute for Nuclear Physics (INFN), grant FLAG.

\section*{Author Contribution}
 Conceptualization F. Mainardi;
methodology, F. Mainardi, A. Mentrelli, J.L. Gonz\'{a}lez-Santander; 
software, A. Mentrelli;
data curation, F. Mainardi, A. Mentrelli, J.L. Gonz\'{a}lez-Santander;
writing---original draft preparation, A. Mentrelli;
writing---review and editing, F. Mainardi, A. Mentrelli, J.L. Gonz\'{a}lez-Santander;
supervision, J.L. Gonz\'{a}lez-Santander.
All authors have read and agreed to the published
version of the manuscript.

\section*{Conflicts of interests}
The authors declare no conflicts of interest.

%%%%%%%%%%%%%%%%%%%%%%%%%%%%%%%%%%%%%%%%%%

%\reftitle{References}


% \bibliography{your_external_BibTeX_file}

%=====================================
% References, variant B: internal bibliography
%=====================================
\newpage
% ACS format
%\isAPAandChicago{}{%

\end{document}